\documentclass[11pt]{article}

\usepackage[dvips]{epsfig} 
\usepackage[dvips]{hyperref}
\usepackage[latin1]{inputenc}
\usepackage{amsmath}
\usepackage{amsfonts}
\usepackage{amssymb}

\def\N{{\mathbb N}}
\def\Hyph#1{\mathbf{H}_{#1}}

\title{A note on strong convergence to common fixed points of nonexpansive mappings 
in Hilbert spaces}

\author{Jean-Philippe Chancelier\thanks{Universit\'e Paris-Est, CERMICS, \'Ecole des Ponts, 6 \& 8 av. B. Pascal, 77455 Marne-la-Vall\'ee, France.}}

\usepackage{a4}
\usepackage[french]{babel}
\usepackage[dvips]{epsfig} 
\usepackage[dvips]{hyperref}
\usepackage[latin1]{inputenc}

\newtheorem{thm}{Theorem}

\newtheorem{lem}[thm]{Lemma}

\newtheorem{prop}[thm]{Proposition}

\newtheorem{defn}[thm]{Definition}

\newtheorem{rem}[thm]{Remark}

\newtheorem{alg}{Algorithm}

\begin{document}

\maketitle 


\begin{abstract}
The aim of this paper is to investigate the links between ${\cal T}_C$-class 
algorithms \cite{bauschke-combettes}, CQ Algorithm \cite{bib18,Solodov-Svaiter-2000} 
and shrinking projection methods \cite{Takahashi-2008}. We show that 
strong convergence of these algorithms are related to coherent 
${\cal T}_C$-class sequences of mapping. Some examples dealing with 
nonexpansive finite set of mappings and nonexpansive semigroups are given. 
They extend some existing theorems in \cite{bauschke-combettes,bib18,Takahashi-2008,saejung-2008}. 
\end{abstract}



\def\sequence#1{(#1_n)_{n\ge 0}}
\def\Fix{\mathop{\mbox{Fix}}}
\def\dom{\mathop{\mbox{dom}}}
\def\VI{\mathop{\normalfont VI}}
\def\defpar{\stackrel{\mbox{\tiny def}}{=}}
\def\argmax{\mathop{\mbox{\rm Argmax}}}
\def\argmin{\mathop{\mbox{\rm Argmin}}}
\def\bbP{{\mathbb P}} 
\def\bbC{{\mathbb C}} 
\def\bbE{{\mathbb E}} 
\def\E{{\cal E}}
\def\F{{\cal F}}
\def\H{{\cal H}}
\def\V{{\cal V}}
\def\W{{\cal W}}
\def\U{{\cal U}}
\def\R{{\mathbb R}}
\def\RB{{\mathbb R}}
\def\N{{\mathbb N}}
\def\M{{\mathbb M}}
\def\S{{\mathbb S}}
\def\bbR{{\mathbb R}} 
\def\bbI{{\mathbb I}} 
\def\bbU{{\mathbb U}} 
\def\Tad{{{\cal T}_{\mbox{\tiny ad}}}}
\def\texte#1{\quad\mbox{#1}\quad}
\def\Proba#1{\bbP\left\{ #1 \right\}} 
\def\Probax#1#2{{\bbP}_{#1}\left\{ #2 \right\}} 
\def\ProbaU#1#2{{\bbP}^{#1} \left\{ #2 \right\}} 
\def\ProbaxU#1#2#3{{\bbP}^{#1}_{#2} \left\{ #3 \right\}} 
\def\valmoy#1{\bbE\left[ #1 \right]}
\def\valmoyDebut#1{\bbE [ #1 } 
\def\valmoyFin#1{ #1 ]} 
\def\valmoyp#1#2{\bbE_{#1}\left[ #2 \right]}
\def\valmoypDebut#1#2{\bbE_{#1} \left[ #2 \right.} 
\def\valmoypFin#1{ \left. #1 \right]} 
\def\valmoypU#1#2#3{\bbE_{#1}^{#2}\left[ #3 \right]}
\def\norminf#1{ {\Vert #1 \Vert}_{\infty}}
\def\norm#1{ {\Vert #1 \Vert}}
\def\psca#1{\left< #1 \right>}
\def\slim{\sigma\mbox{-}\lim}

\newenvironment{myproof}{{\small{\it Proof:}}}{\hfill$\Box$\normalsize
\\\smallskip}

\section{Introduction} 

Let $C$ be a closed convex subset of a Hilbert space $\H$. A mapping $T$ of $C$
into itself is called nonexpansive if
\begin{equation*}
  \norm{Tx -Ty} \le \norm{x-y} \quad \mbox{for all} x,y \in C\,.
\end{equation*}
We denote by $\Fix(T)$ the set of fixed points of $T$. That is
\begin{equation}
  \Fix(T) \defpar  \left\{x \in C\, :\,  T x = x \right\}\,.
\end{equation}

There are many iterative methods for approximation of common fixed points of a
family of nonexpansive mapping in a Hilbert space.  In Section \ref{secacedo} we
recall the CQ Algorithm \cite{bib18,Solodov-Svaiter-2000} (Algorithm
\ref{firstalg}) associated to a sequence of mappings $\sequence{T}$ of $C$ into
itself. The CQ Algorithm when applied to a sequence of mappings of $\H$ into
itself is the same as a Haugazeau method \cite{Haugazeau} studied in
\cite[Algorithm 3.1]{bauschke-combettes} and applied to ${\cal T}$-class
mappings.

We straighforwardly generalize, in Section \ref{secacedo}, the ${\cal T}$-class
to take into account mappings of $C$ into itself. We denote this new class by
the ${\cal T}_C$-class.  Using this extension, the CQ Algorithm (Algorithm
\ref{firstalg}) coincides with the Haugazeau method (Algorithm \ref{bcalg}) and
a strong convergence theorem can be obtained by following results from
\cite{bauschke-combettes}.  Note that the convergence theorem is obtained for
${\cal T}_C$-class sequences which are coherent (Definition \ref{coherentprop}).

In \cite{Takahashi-2008} another algorithm called the shrinking projection
method is also studied. One of our aims in this article is to prove that,
rephrased in the context of ${\cal T}_C$-class algorithm, the convergence
results of this new algorithm (Algorithm \ref{newalg}) is also related to
coherent sequences of ${\cal T}_C$-class mappings. We give in Theorem
\ref{maint_th} a strong convergence result of Algorithm \ref{newalg} for ${\cal
  T}_C$-class coherent sequence of mappings. Section \ref{proof} is devoted
to the proof.  The strong convergence of Algorithm \ref{newalg} is also proved
in \cite[Theorem 3.3]{Takahashi-2008} for sequence of nonexpansive mappings
satisfying the NST-condition(I) (Definition \ref{defnst}).  It is easy to prove
that if $R$ is a nonexpansive mapping of $C$ into itself then $T=(R+Id)/2$ belongs
to the ${\cal T}_C$-class and that a sequence of nonexpansive mappings satisfying
the NST-condition(I) is coherent. Thus Theorem \ref{maint_th} extends \cite[Theorem
3.3 and Theorem 3.4]{Takahashi-2008}.

In Section \ref{examples} we show that specific sequences of mappings are 
coherent. Combined with  Theorem \ref{maint_th} it can be considered as an 
extension to some existing theorems in \cite{bib18,Takahashi-2008,saejung-2008}. 

\section{The ${\cal T}_C$-class iterative algorithms,  CQ algorithm and the shrinking projection method}

\label{secacedo}

We first recall here the ${\cal T}$-class iterative algorithms as defined 
by H. Bauschke and P. L. Combettes \cite{bauschke-combettes}.

For $(x,y)\in \H^2$ and $S$ a subset of $\H$, we define the mapping $H_S$ as
follows:
\begin{equation}
  H_S(x,y) \defpar \left\{ z \in S \quad \vert \quad 
    \psca{ z -y , x - y} \le 0  \right\}\,.
\end{equation}
We also define the mapping $H$ by $H \defpar H_{\H}$. Note that $H_S(x,x)=S$ and
for $x \ne y$, $H(x,y)$ is a closed affine half space. For a nonempty closed
convex $C$, we denote by $Q_C(x,y,z)$ the projection, when it exists, of $x$ onto
$H_C(x,y)\cap H_C(y,z)$ and $Q$ the projection when $C=\H$, that is $Q \defpar
Q_{\H}$. As an intersection of two closed affine half spaces and a closed convex,
$H_C(x,y)\cap H_C(y,z)$ is a possibly empty closed convex.  

It is easy to check, from the definition of $H$, that $y$ is the projection 
of $x$ onto $H(x,y)$ and we therefore have $Q(x,x,y)= P_{H(x,y)} x = y$. 
Where $P_C$ is the metric projection from $\H$ onto $C$.
Moreover, if $y\in C$ then we also have that $y$ is the 
projection of $x$ onto $H_C(x,y)$ which gives $Q_C(x,x,y)=y$.

The algorithm studied in \cite{bauschke-combettes} is the following

\begin{alg}
  \label{bcalg}
  Given $x_0 \in C$ and a sequence $\sequence{T}$ of mappings $T_n : C \to \H$,
  we consider the sequence $\sequence{x}$ generated by the following algorithm:
\begin{equation}
  x_{n+1}= Q_C(x_0,x_n,T_{n} x_n) \nonumber 
\end{equation}
\end{alg}

A very similar algorithm exists under the name of CQ algorithm \cite{bib18,Solodov-Svaiter-2000}:

\begin{alg}
\label{firstalg}
Given $x_0 \in C$, we consider the sequence $\sequence{x}$ generated by the following algorithm:
\begin{equation}
  \begin{cases}
    y_n = R_n x_n, \nonumber \\ 
    C_n \defpar \left\{ z \in C \,| \, \norm{y_n - z} \le \norm{x_n - z}\right\} \,, \nonumber \\ 
    D_n \defpar \left\{ z \in C \,|\, \psca{x_n - z, x_0 - x_n} \ge 0 \right\}, \nonumber \\ 
    x_{n+1} = P_{(C_n \cap D_n)} x_0.\nonumber
\end{cases}
\end{equation}
\end{alg}

The link between the two algorithms is described by the following lemma.

\begin{lem}The sequence generated by Algorithm \ref{firstalg} coincides with the
  sequence given by $x_{n+1}= Q_C(x_0,x_n,T_{n} x_n)$ with $T_n \defpar (R_n
  +Id)/2$.
\end{lem}

\begin{myproof} Following \cite{bauschke-combettes},  the proof easily follows from the equality 
\begin{equation}
  4 \psca{ z -T x, x - Tx}= \norm{Rx - z }^2 - \norm{x-z}^2 \,. \nonumber 
\end{equation}
\end{myproof}

The convergence of Algorithm \ref{bcalg} and therefore of Algorithm
\ref{firstalg} when $C=\H$ is studied in \cite{bauschke-combettes}. It relies
on two requested properties of the sequence $\sequence{T}$. First, the
 sequence $\sequence{T}$ must belong the ${\cal T}$-class which means that for all $n\in
\N$ we must have $T_n \in {\cal T}$ where ${\cal T}$ is defined as follows:

\begin{defn} A mapping $T: C \mapsto \H$ belongs to the 
  ${\cal T}_C$-class if it is an element of the set ${\cal T}_C$:
\begin{equation}
  {\cal T}_C \defpar \left\{ T: C \mapsto C \, | \, \dom(T) = C \quad\mbox{and} \quad (\forall x \in C) 
  \Fix(T) \subset H(x,Tx)\right\} \nonumber\,.
\end{equation}
\end{defn}
When $C=\H$, we use the notation ${\cal T}= {\cal T}_{\H}$.
Second, the sequence $\sequence{T}$ must be coherent as defined below.

\begin{defn}\cite{bauschke-combettes} A sequence $\sequence{T}$ such that $T_n \in {\cal T}_C$ is {\emph coherent} if for every bounded sequence $\{z_n\}_{n\ge 0} \in C$ the
 following holds: 
\begin{eqnarray}
  \left\{ 
  \begin{array}{l} 
    \sum_{n\ge 0} \norm{z_{n+1} -z_n}^2 < \infty \\
    \sum_{n\ge 0} \norm{z_{n} - T_n z_n}^2 < \infty 
  \end{array} 
  \right.
  \Rightarrow {\cal M}(z_n)_{n\ge 0} \subset \bigcap_{n \ge 0} \Fix(T_n)
  \label{coherentprop}
\end{eqnarray} 
where ${\cal M}(z_n)_{n\ge 0}$ is the set of weak cluster points of the sequence $\sequence{z}$.
\end{defn}

\begin{thm}\cite[Theorem 4.2]{bauschke-combettes} Suppose that $C=\H$ and the ${\cal T}_C$-class sequence $\sequence{T}$ is coherent. Then, for an arbitrary 
  orbit of Algorithm \ref{bcalg}, exactly one of the following alternatives
  holds:
\begin{enumerate}
  \item[(a)] $F\ne \emptyset$ and $x_n \to_{n} P_F x_0$;
  \item[(b)] $F=  \emptyset$ and $x_n \to_{n} +\infty $;
  \item[(c)] $F=  \emptyset$ and the algorithm terminates,
\end{enumerate}
where the set $F$ is defined by $F\defpar \bigcap_{n \ge 0} \Fix(T_n)$.
\label{bcth}
\end{thm}

\begin{rem}
In the previous proof, it is supposed that $C=\H$. 
If $C$ is a nonempty closed convex subset of
$\H$, Theorem \ref{bcth} $(a)$ remains valid. 
\end{rem}

In \cite{Takahashi-2008} another iterative algorithm called the 
\emph{shrinking projection method} is studied. Using our notation
it can be rephrased as follows:

\begin{alg}
\label{newalg}
Given $x_0 \in C$ and $C_0 \defpar C$, we consider the sequence $\sequence{x}$ 
(when it exists) generated by the following algorithm:
\begin{equation}
\begin{cases}
  C_{n+1} \defpar C_n \cap H(x_n,T_n x_n) \quad\mbox{with}\quad 
  T_n \defpar (R_n +Id)/2\,, \nonumber \\ 
  x_{n+1} = P_{C_{n+1}} x_0.\nonumber
\end{cases}
\end{equation}
\end{alg}

The previous algorithm is stopped once $C_n = \emptyset$. 
One of the results of this paper is the proof that the convergence of Algorithm \ref{newalg}
is governed by the same rules as for the convergence of Algorithm \ref{bcalg}. 

\begin{thm} Suppose that the ${\cal T}_C$-class sequence $(T_n)_{n\in \N}$ is
  coherent and let \[ F \defpar \bigcap_{n \in \N} \Fix(T_n)\,.\]
 Then, if $F\ne \emptyset$ the sequence $\sequence{x}$ produced by Algorithm \ref{newalg} and
 Algorithm \ref{bcalg} converges to $P_F x_0$.
 \label{maint_th}
\end{thm}

\begin{myproof} As pointed out in the introduction the case of  Algorithm \ref{bcalg} when 
  $C=\H$ is proved in Theorem \ref{bcth}. The extension to the case of a closed nonempty subset 
  $C$ of $\H$ is straightforward and we will not give an explicit proof. The proof 
  of the case of  Algorithm \ref{newalg} is postponed to Section \ref{proof}.
\end{myproof}

\begin{rem} The first condition for the convergence is the fact that the 
sequence $\sequence{T}$ must belong to the ${\cal T}_C$-class. Note 
that by \cite[Proposition 2.3]{bauschke-combettes} $T \in {\cal T}$ iff 
the mapping $2T-Id$ is quasi nonexpansive and $\dom(T)=\H$. The equivalence 
remains true for ${\cal T}_C$-class if $\dom(T)=\H$ is replaced by 
$\dom(T)=C$. 

Thus, if $T_n \defpar (R_n +Id)/2$, a necessary and sufficient condition for the
sequence $\sequence{T}$ to belong to the ${\cal T}_C$-class is that the sequence
$\sequence{R}$ is a sequence of quasi nonexpansive mappings.
\end{rem}

\begin{rem} Moreover, it is a well known fact \cite[Theorem 12.1]{bib5} that 
$2T-Id$ is nonexpansive iff $T$ is firmly nonexpansive. Thus, a sufficient 
condition for the mapping $T$ to belong to the ${\cal T}_C$-class is that $T$ is a {\em firmly 
nonexpansive} mapping, {\em i.e}: 
\begin{equation}
  \norm{T x - Ty }^2 \le \psca{x-y,T x - Ty} \quad \forall (x,y)\in C^2
\end{equation}
or equivalently 
\begin{equation}
  \norm{T x - Ty }^2 \le \norm{x-y}^2 - \norm{(T-Id)x -(T-Id)y}^2 \quad \forall (x,y)\in C^2 \,.
\end{equation}
\end{rem}


We recall here the definition of the NST-condition (I) \cite{Nakajo-Shimoji-Takahashi}.
Let $\sequence{T}$ and ${\cal F}$ be two families of nonexpansive mappings of $C$
into itself such that
\[
\emptyset \ne \Fix({\cal F}) \defpar  \bigcap_{n \in \N} \Fix(T_n) \,,
\] 
where $\Fix({\cal F})$ is the set of all common fixed points of mappings from the family
${\cal F}$. 
\begin{defn} The sequence $\sequence{T}$ of mappings is said to satisfy the 
NST-condition (I) with ${\cal F}$ if, for each bounded sequence $\sequence{z} \subset C$, 
we have that $\lim_{n\mapsto \infty} \norm{z_n -T_n z_n}= 0$ implies that 
$\lim_{n\mapsto \infty} \norm{z_n -T z_n}= 0$ for all $T \in {\cal F}$. 
\label{defnst}
\end{defn}

\begin{rem}Suppose that ${\cal F}$ is a family of nonexpansive mappings. 
It is easy to see that a sequence $\sequence{T}$ of mappings 
satisfying a NST-condition (I) with ${\cal F}$ is coherent. Indeed, from a demi-closed principle or using \cite[Lemma 3.1]{Takahashi-2008} if $\norm{x_n -T x_n} \mapsto 0$ for all $T \in {\cal T}$ then ${\cal M}(x_n)_{n\ge 0} 
\subset \Fix\left( \{ T \}_{T \in {\cal T}} \right)$. 
\end{rem}

\section{Coherent sequences of mappings}

\label{examples}

We consider here Algorithms \ref{bcalg} and \ref{newalg} for 
a sequence of mappings $\sequence{R}$ built by $N$ level iterations. 
Our aim is to give conditions under which the 
sequence $\sequence{R}$ or equivalently $\sequence{T} \defpar (R_n +Id)/2$ 
is coherent\footnote{By \cite[Proposition 4.5]{bauschke-combettes} if 
$\sequence{T} \in {\cal T}$ and $T'_n \defpar Id + \lambda_n(T_n -Id)$ with 
$\lambda_n \in [\delta ,1]$ and $\delta \in ]0,1]$. Then $\sequence{T}$ 
is coherent iff $\sequence{T'}$ is coherent.} and apply Theorem \ref{maint_th}
to get convergence results. 

Let $N \ge 1$ and $\sequence{T^{(j)}} : C \to \H$ for $1 \le j \le N$ be a finite 
set of sequences of nonexpansive mappings. Given also a family of 
sequences of real parameters $\sequence{\alpha^{(j)}}$ for $1 \le j \le N$,
we define new sequences $\sequence{\Gamma^{(j)}} : C \to \H$ by the 
recursive equations:

\begin{equation}
  \Gamma^{(j)}_n x \defpar \alpha^{(j)}_n x + (1 - \alpha^{(j)}_n) T^{(j)}_n \Gamma^{(j+1)}_n x 
  \quad\mbox{and}\quad \Gamma_n^{(N+1)} x \defpar x
  \label{gammandef}
\end{equation}

$\Hyph{\alpha}$: We will assume that the sequences of real parameters $\sequence{\alpha^{(j)}}$ satisfy the following condition:
for $2 \le j \le N$ and for all $n\in \N$ we have $\alpha^{(j)}_n \in (a,b)$ with $0 < a < b < 1$ 
and $\alpha^{(1)}_n \in [0,b)$. 

Using the sequence of mappings $R_n \defpar \Gamma_n^{(1)}$ in Algorithms \ref{bcalg} 
and \ref{newalg} gives $N$ level algorithms. We will consider the following specific 
examples:

\begin{description}
  \item[$\Hyph{1}$]  Each sequence $\sequence{T^{(j)}}$ is constant, 
    {\emph i.e} $T^{(j)}_n = T^{(j)}$ for $  1 \le j \le N$ and 
    $F \defpar \Fix\left( \left\{T^{(j)}, 1 \le j \le N \right\} \right)$ is nonempty.
  \item[$\Hyph{2}$] The $\sequence{T^{(j)}}$ sequences for $1 \le j \le N$ are
    given by $T^{(j)}_n = T^{(j)}(t_n)$, where $\left\{ T^{(j)}(t): t \ge
      0\right\}$ is a finite set of given semigrougs and $\sequence{t}$ is a
    sequence of real numbers such that $\liminf_n t_n =0$, $\limsup_n t_n >0$
    and $\lim_n (t_{n+1}-t_n)=0$.  We assume that $F \defpar \Fix \left( \left\{
        T^{(j)}(t), 1 \le j \le N , t \ge 0 \right\} \right)$ is nonempty.
  \item[$\Hyph{3}$] The $\sequence{T^{(j)}}$ sequences for $1 \le j \le N$ are given by 
    \begin{equation}
      T^{(j)}_n x = \frac{1}{t_n} \int_0^{t_n} T^{(j)}(s) x ds \,,
      \label{eqdeftnh3}
    \end{equation}
      where $\left\{ T^{(j)}(t): t \ge 0\right \}$ is a finite set
      of given semigrougs and $\sequence{t}$ is a positive divergent sequence 
      of real numbers. We assume that $F \defpar \Fix\left( \left\{T^{(j)}(t), 1 \le j \le N
         , t \ge 0 \right\} \right)$ is nonempty.
\end{description}

\begin{thm} Given a finite set of $N$ nonexpansive sequences $\sequence{T^{(j)}}$ satisfying
  $\Hyph{1}$, $\Hyph{2}$, or $\Hyph{3}$. 
  The sequence $\sequence{x}$ produced by Algorithm \ref{bcalg} and
  Algorithm \ref{newalg} with $R_n \defpar \Gamma_n^{(1)}$ and 
  $\sequence{T}\defpar (R_n +Id)/2$ converges to $P_F x_0$.   
  The mappings $\Gamma_n^{(j)}$ being defined by equation \eqref{gammandef} with 
  parameters $\alpha_n^{(j)}$ satisfying $\Hyph{\alpha}$.
\end{thm} 

\begin{myproof} The proof is obtained by showing that the sequence of 
  mappings $\sequence{T}$ is coherent in each given case and 
  by applying Theorem \ref{maint_th} to conclude. 
  The coherence is proved in the sequel in Proposition \ref{prop-coherence} 
  for the case $\Hyph{1}$, in Proposition \ref{prop-coherence-semig} for the case $\Hyph{2}$ and in 
  Proposition \ref{prop-coherence-semig-last} for the case $\Hyph{3}$.
\end{myproof}

We start here by a set of lemmata which are common to all cases. 

\begin{lem} Let $T$ be a $F$-quasi nonexpansive mapping and for $\beta \in (0,1)$ 
the mapping $T_{\beta} \defpar \beta Id + (1-\beta)T$. For $p\in F$ and all $x\in H$ 
we have~:
\begin{equation}
  \beta(1-\beta)\norm{ x - T x}^2 \le 2 ( \norm{x-p}-\norm{T_{\beta} x -p})\norm{x-p}
\end{equation}
\label{lem:un}
\end{lem}

\begin{myproof} For $p\in F$ and all $x\in H$ we have~:
\begin{eqnarray}
  \norm{T_{\beta} x -p}^2 &= & \norm{\beta(x-p)+(1-\beta) (Tx-p)}^2 \nonumber \\
  &= & \beta\norm{x-p}^2 + (1-\beta)\norm{Tx-p}^2 - \beta(1-\beta) 
  \norm{Tx -x}^2 \nonumber \\ 
  &\le & \norm{x-p}^2 - \beta(1-\beta)\norm{Tx -x}^2 \nonumber \,.
\end{eqnarray} 
We thus obtain
\begin{eqnarray}
  \beta(1-\beta)\norm{Tx -x}^2 &\le & ( \norm{x-p}-\norm{T_{\beta} x -p})( \norm{x-p}+\norm{T_{\beta} x -p}) \nonumber \\
  & \le & 2 ( \norm{x-p}-\norm{T_{\beta} x -p})\norm{x-p}\,. \nonumber
\end{eqnarray} 
\end{myproof}

\begin{lem} Let $T$ a $F$-quasi nonexpansive mapping. For $\beta \in (0,1)$ we define the mapping  $T_{\beta} \defpar \beta Id + (1-\beta)T$. For $p\in F$, all $x\in H$ and $S$ a $F$-quasi nonexpansive mapping, we have~:
\begin{equation}
  \beta(1-\beta)\norm{ x - T x}^2 \le 2 \norm{x- S T_{\beta} x} \norm{x-p}
  \label{eq:lemdeuxun}.
\end{equation}
If moreover $S$ is nonexpansive we also have~:
\begin{equation}
  \norm{ x - S x} \le \norm{x- S T_{\beta} x} + \norm{Tx-x}
  \label{eq:lemdeuxdeux}.
\end{equation}
\label{lemdeux}
\end{lem}

\begin{myproof} For $p\in F$ and all $x\in H$ we have~:
\begin{eqnarray}
  \norm{x -p} &\le & \norm{ x- ST_\beta x} + \norm{ST_\beta x -p} \nonumber \\
  &\le  & \norm{ x- ST_\beta x} + \norm{T_\beta x -p} \nonumber \,.
\end{eqnarray}
We thus have $ \norm{x -p} - \norm{T_\beta x -p} \le \norm{ x- ST_\beta x}$ which 
combined with Lemma \ref{lem:un} gives equation \eqref{eq:lemdeuxun}.

Now if $S$ is nonexpansive, 
\begin{eqnarray}
  \norm{x - Sx} &\le & \norm{ x - ST_\beta x} + \norm{ST_\beta x - Sx } 
  \le  \norm{ x- ST_\beta x} + \norm{T_\beta x -x} \nonumber \\
  & \le & \norm{ x- ST_\beta x} + (1-\beta) \norm{T x -x} 
  \le  \norm{ x- ST_\beta x} + \norm{T x -x} \nonumber\,.
\end{eqnarray}
\end{myproof}

\begin{lem} Suppose that $F \defpar \bigcap_{ \{ n\in \N; 1 \le j \le N\}}
  \Fix(T^{(j)}_n)$ is not empty \and suppose that for a bounded sequence
  $\sequence{x}$ and a fixed value of $j$ we have\\ $\norm{ x_n - T^{(j)}_n
    \Gamma_n^{(j+1)} x_n} \to 0$. Moreover, suppose that for $2 \le j \le N$ and
  all $n\in \N$ we have $\alpha^{(j)}_n \in (a,b)$ with $0 < a < b < 1$.  Then
  for all $k \ge j$ we have $\norm{ x_n - T^{(k)}_n x_n} \to 0$.
\label{lem-coherence}
\end{lem}

\begin{myproof} Note first that the sequences $(T^{(j)})_{ 1 \le j \le N}$ and
  $(\Gamma^{(j)})_{ 1 \le j \le N+1}$ are composed of nonexpansive mappings.
  Indeed the composition of nonexpansive mappings is nonexpansive and for
  $\beta \in (0,1)$ $\beta Id + (1- \beta) S$ is nonexpansive when $S$ is
  nonexpansive. The sequences are also $F$-quasi nonexpansive since it is
  straightforward that $F \subset \Fix(\Gamma^{(j)}_n)$ for all $j \in[1,N]$ and 
  $n\in \N$ and if $S$ is nonexpansive it is also $\Fix(S)$-quasi nonexpansive. 

  The proof then follows by backward induction on $j$. Assume that the result is
  true for $j+1$ then we will prove that it is true for $j$. Using the
  definition of $\Gamma_n^{(j+1)}$ and using equation \eqref{eq:lemdeuxun} 
  for $p\in F$, $S=T^{(j)}_n$, $T=T^{(j+1)}_n \Gamma_n^{(j+2)}$ 
  and $\beta =\alpha_n^{(j+1)}$ (we thus have $T_{\beta}=\Gamma_n^{(j+1)}$)
  we obtain~:
  \begin{equation}
    \alpha_n^{(j+1)} (1-\alpha_n^{(j+1)}) \norm{ x_n - T^{(j+1)}_n \Gamma_n^{(j+2)} x_n}^2 \le 2 \norm{ x_n - T^{(j)}_n \Gamma_n^{(j+1)} x_n} \norm{x_n-p}
  \end{equation}
  We thus obtain that $\norm{ x_n - T^{(j+1)}_n \Gamma_n^{(j+2)} x_n} \to 0$ and by
  induction hypothesis we obtain $\norm{ x_n - T^{(k)}_n x_n} \to 0$ for $k \ge j+1$.
  Now using equation \eqref{eq:lemdeuxdeux} with 
  $S\defpar T^{(j)}_n$, $T \defpar T^{(j+1)}_n\Gamma_n^{(j+2)}$ and $\beta =\alpha_n^{(j+1)}$
  we get:
\begin{equation}
  \norm{ x_n - T^{(j)}_n x_n} \le \norm{x_n - T^{(j)}_n \Gamma_n^{(j+1)} x_n} +
  \norm{T^{(j+1)}_n \Gamma_n^{(j+2)}x_n-x_n}
\end{equation}
and the result follows for $j$.
\end{myproof}

\subsection{The case $\Hyph{1}$}

\begin{prop} In the case $\Hyph{1}$, the sequence $\sequence{R}$, defined by $R_n \defpar \Gamma^{(1)}_n$ with parameters satisfying $\Hyph{\alpha}$, satisfy the NST-condition(I) 
with \\ ${\cal F}\defpar \Fix{ \{ T^{(j)}}_{1 \le j \le N}\}$ and the sequence 
$T_n = (R_n +Id)/2$ is a ${\cal T}_C$-class and coherent sequence. 
\label{prop-coherence}
\end{prop}

\begin{myproof} We have $\norm{x_n - R_n x_n}= \norm{x_n - T_n^{(1)} \Gamma_n^{(2)}x_n}(1 -\alpha^{(1)}_n)$. Thus, if for each bounded sequence $\sequence{x}$ $\norm{x_n - R_n x_n}\mapsto 0$ we also 
have $\norm{x_n - T_n^{(1)} \Gamma_n^{(2)}x_n} \mapsto 0$ since $(1 -\alpha^{(1)}_n)$ is bounded 
from zero. Using Lemma \ref{lem-coherence} we have $\norm{ x_n - T^{(j)} x_n} \mapsto 0$ 
for $1 \le j \le N$ which gives use the NST-condition(I) with ${\cal F}$. 
Now we consider the sequence $\sequence{T}$. The sequence belongs to the ${\cal T}_C$-class 
since $2T_n-Id = R_n$ is nonexpansive and thus quasi nonexpansive. Now if 
$\norm{x_n - T_n x_n} \mapsto 0$ we also have $\norm{x_n - R_n x_n} \mapsto 0$ and thus 
using the NST-condition(I) we have $\norm{ x_n - T^{(j)} x_n} \mapsto 0$ for $1 \le j \le N$. 
Since the $T^{(j)}$ are nonexpansive they are also demi-closed \cite[Lemma 4]{browder-1967} and 
thus we must have ${\cal M}(x_n)_{n\ge 0} \subset \Fix( \{ T^{(j)}, 1 \le j \le N\} ) = \Fix( \{ T_n \}_{n \in \N})$. 
The sequence $\sequence{T}$ is thus in the ${\cal T}_C$-class and coherent.
\end{myproof}

\begin{rem} For $N=1$ we recover \cite[Theorem 1.1]{Takahashi-2008} and 
  \cite[Theorem 4.1]{Takahashi-2008}. 
\end{rem}

\subsection{The case $\Hyph{2}$}

Let $\left\{T(t) : t \ge 0\right\}$ be a family of mappings from a subset $C$ of $\H$
    into itself. We call it a nonexpansive semigroup on $C$ 
    if the following conditions are satisfied:
\begin{itemize}
\item[(i)] $T(0)x = x$ for all $x \in C$;
\item[(ii)] $T(s + t) = T(s)T(t)$ for all $s, t \ge 0$;
\item[(iii)] for each $x \in C$ the mapping $t \mapsto T(t)x$ is continuous;
\item[(iv)] $\norm{T(t)x - T(t)y} \le \norm{x - y}$ for all $x, y \in C$ and 
  $t \ge 0$.
\end{itemize}

\begin{prop}  In the case $\Hyph{2}$, the sequence $\sequence{R}$, defined by $R_n \defpar \Gamma^{(1)}_n$ with parameters satisfying $\Hyph{\alpha}$, satisfy the NST-condition(I) 
  with \\ ${\cal F}\defpar \Fix{ \{
    T^{(j)}(t)}_{1 \le j \le N, t \ge 0}\}$ and the sequence $T_n = (R_n +Id)/2$ is a ${\cal
    T}_C$-class and coherent sequence.
\label{prop-coherence-semig}
\end{prop}

\begin{myproof} As in the proof of Proposition \ref{prop-coherence} we obtain
  that for each bounded sequence $\sequence{x}$ such that $\norm{x_n - R_n
    x_n}\mapsto 0$ we also have $\norm{x_n - T^{(j)}(t_n) x_n} \mapsto 0$ for $1
  \le j \le N$.  Now it is easy to prove that the weak cluster points of the
  sequence $\sequence{x}$ are in $F$. The proof for each fixed $j$ is the same
  as in \cite[Theorem 2.2, page 6]{saejung-2008}. We thus obtain the coherence
  of the sequence $\sequence{T}$.
\end{myproof}

\begin{rem}
For $N=1$ we recover \cite[Theorem 2.1]{saejung-2008} for Algorithm \ref{newalg} and \cite[Theorem 2.2]{saejung-2008} for Algorithm \ref{bcalg}.
\end{rem}

\subsection{The case $\Hyph{3}$} 

\begin{prop}  In the case $\Hyph{3}$, the sequence $\sequence{R}$, defined by $R_n \defpar \Gamma^{(1)}_n$ with parameters satisfying $\Hyph{\alpha}$, satisfy the NST-condition(I) 
  with \\ ${\cal F}\defpar \Fix{ \{
    T^{(j)}(t)}_{1 \le j \le N, t \ge 0}\}$ and the sequence $T_n = (R_n +Id)/2$ is a ${\cal
    T}_C$-class and coherent sequence.
\label{prop-coherence-semig-last}
\end{prop}

\begin{myproof} As in the proof of Proposition \ref{prop-coherence} we obtain that 
for each bounded sequence $\sequence{x}$ such that $\norm{x_n - R_n x_n}\mapsto 0$ 
we also have $\norm{x_n - T^{(j)}(t_n) x_n} \mapsto 0$ for $1 \le j \le N$. 
Now it is easy to prove that the weak cluster points of the sequence $\sequence{x}$ 
are in $F$. The proof for each fixed $j$ is the same as in \cite[Theorem 4.1]{bib18}. 
For each fixed $j$, it is a consequence of the inequality \cite[Equation (8)]{bib18}:
\begin{equation} 
  \norm{T^{(j)}(s) x_n -x_n} \le 2 \norm{ T^{(j)}_n x_n -x_n} 
  + \norm{ T(s) ( T^{(j)}_n x_n) - T^{(j)}_n x_n} 
\end{equation}
for every $0 \le s < + \infty$ and $n \in \N$ with $T^{(j)}_n$ and the fact that 
the right hand side of the above inequality goes to zero as $n$ goes to infinity for a bounded 
sequence $\sequence{x}$ using \cite[Lemma 2.1]{bib18}.
We thus obtain the coherence of the sequence $\sequence{T}$. 
\end{myproof}

\begin{rem}
For $N=1$ we recover \cite[Theorem 4.1]{bib18} for Algorithm \ref{bcalg} and \cite[Theorem 4.4]{Takahashi-2008} for Algorithm \ref{newalg}.
\end{rem}

\section{Proof of Theorem \ref{maint_th}}

\label{proof}

We prove here the strong convergence of Algorithm \ref{newalg} for a ${\cal T}_C$-class sequence 
of coherent mappings. The proof follows the same steps as the proof of the convergence of  Algorithm 
\ref{bcalg} in \cite{bauschke-combettes}, we therefore give references to the original propositions. 

The proof results from the next proposition and theorem in the following way.
Let $\sequence{x}$ be an arbitrary orbit of Algorithm \ref{newalg} and let 
$F \defpar \Fix( \{ T_n\}_{n \in \N} )$. If $F\ne \emptyset$, then by Proposition 
\ref{prop:tq} $(iv)$ the sequence is defined. By Theorem \ref{thm35} $(ii)$ the sequence 
is bounded. Thus $(v)$ is fulfilled and by the coherence property we have 
$ {\cal M}(x_n)_{n\ge 0} \subset F$. Then, by Theorem \ref{thm35} $(iv)$, the sequence 
strongly converges to $P_F (x_0)$. 

\begin{prop} \cite[Proposition 3.4]{bauschke-combettes} Let $\sequence{x}$ be an arbitrary 
orbit of Algorithm \ref{newalg}. Then:
\begin{itemize}
\item[$(i)$] If $x_{n+1}$ is defined then $\norm{x_0-x_n} \le \norm{x_0 -x_{n+1}}$.
\item[$(ii)$] If $x_{n}$ is defined then $x_0=x_n \iff x_n =x_{n-1}=\cdots = x_0 \iff 
  x_0 \in \bigcup_{k=0}^{n-1} \Fix(T_k)$.
\item[$(iii)$] If $\sequence{x}$ is defined then $ (\norm{x_0-x_n})_{n\in \N}$ is increasing.
\item[$(iv)$] $\sequence{x}$ is defined if $F \defpar \Fix( \{ T_n\}_{n \in \N} )\ne \emptyset$.
\end{itemize}
\label{prop:tq}
\end{prop}

\begin{myproof} $(i)$: If $x_{n+1}$ is defined we have $x_{n+1} = P_{C_{n+1}}
  x_0$ and thus $x_{n+1} \in C_{n+1} \subset C_n$ and since $x_n = P_{C_{n}}
  x_0$ we have $\norm{x_0-x_n} \le \norm{x_0 -x_{n+1}}$.
  $(ii)$: The fist equivalence follows from $(i)$. The second one is proved 
  by induction. Note first that $H$ is such that $y = P_{H(x,y)} x$. Now for $y \in C$, 
  we obtain also that $y = P_{C \cap H(x,y)} x$. for $n=1$, we have 
  $x_1 = P_{C \cap H(x_0,T_0 x_0)} x_0 = T_0 x_0$ and thus $x_1=x_0 \iff x_0 \in \Fix(T_0)$. 
  Now assume that the equivalence if fulfilled for $n$. We have 
  \[
   x_{n+1} =x_{n}=\cdots = x_0 \iff 
   \begin{cases}
     x_0 &\in \cup_{k=0}^{n-1} \Fix(T_k)\, \\
     x_0& = x_{n+1} = P_{ C \cap \bigcap_{k=0}^n H(x_k,T_k x_k)} \\
     & =P_{C \cap H(x_0,T_n x_0)} = T_n x_0    \,.
   \end{cases}
  \]
  $(iii)$ follows from $(i)$. $(iv)$: The algorithm is defined if $C_n \ne \emptyset$ for 
  all $n \in \N$. Thus it is enough to prove that $ C \cap \left(\bigcap_{n \in \N} H(x_n,T_n x_n)\right) \ne \emptyset$. 
  By definition of the ${\cal T}_C$ class we have $\Fix(T_n) \subset C \cap H(x_n,T_n x_n)$ and the 
  result follows. 
\end{myproof}

\begin{thm} (\cite[Theorem 3.5]{bauschke-combettes}) Let $\sequence{x}$ be an arbitrary 
orbit of Algorithm \ref{newalg} and let $F \defpar  \bigcap_{n \in \N} \Fix(T_n)$. Then 
\begin{itemize}
\item[$(i)$] If $\sequence{x}$ is defined then: 
  $\sequence{x}$ is bounded $\iff$ $ (\norm{x_0-x_n})_{n\in \N}$ converges. 
\item[$(ii)$] If $F\ne \emptyset$, then  $\sequence{x}$ is bounded and $(\forall n \in \N) 
  x_n \in F \iff x_n =P_F (x_0)$.
\item[$(iii)$]  If $F\ne \emptyset$, then $ (\norm{x_0-x_n})_{n\in \N}$ converges and \\
  $\lim_n \norm{x_0 -x_n} \le \norm{ x_0 -P_F x_0}$.
\item[$(iv)$]  If $F\ne \emptyset$, then:  $ \lim_n x_n = P_F(x_0) \iff {\cal M}(x_n)_{n \in \N} \subset F$.
\item[$(v)$] If  $\sequence{x}$ is defined and bounded then 
  $\sum_{n \ge 0} \norm{x_{n+1}-x_n}^2 < + \infty$ and 
  $\sum_{n \ge 0} \norm{x_{n}- T_n x_n}^2 < + \infty$. 
\end{itemize}
\label{thm35}
\end{thm}

\begin{myproof} $(i)$ follows from Proposition \ref{prop:tq} $(i)$. $(ii)$:  If $F\ne \emptyset$ 
then by Proposition \ref{prop:tq} $(iv)$ the sequence is defined. We have $F \subset C 
\cap \left( \bigcap_{n\in \N} H(x_n,T_n x_n) \right)$ and thus $F \subset C_{n}$. Now, from $P_F (x_0) \in C_{n}$  
and $x_{n}=P_{ C_n} x_0$ we obtain $\norm{x_n -x_0} \le \norm{x_0 -P_F(x_0)}$ and $(ii)$ follows. 
$(iii)$ follows from $(i)$, $(ii)$ and the previous inequality. $(iv)$: The forward implication
is trivial. For the reverse implication, the proof exactly follows $(iv)$ of \cite[Theorem 3.5]{bauschke-combettes} since it does not involve $C$. 
$(v)$: From $x_n = P_{C_n} x_0$ and $x_{n+1}\in C_n$ we obtain:
\[
\psca{ x_0 -x_n, x_n-x_{n+1} } \ge 0\,.
\]
We thus have~: 
\begin{eqnarray}
  \norm{x_{n+1} -x_n}^2 &\le & \norm{x_{n+1} -x_n}^2 + 2\psca{x_{n+1} -x_n, x_n -x_0} \nonumber \\
  & \le & \norm{x_{0} -x_{n+1}}^2 -  \norm{x_{0} -x_{n}}^2 \,.
\end{eqnarray}
Hence $\sum_{n \ge 0} \norm{x_{n+1} -x_n}^2 \le \sup_{n\in \N} \norm{x_{0} -x_{n}}^2 < + \infty$ since 
$\sequence{x}$ is bounded. For all $n \in \N$ we have $x_{n+1} \in H(x_n,T_n x_n)$, which implies,
\begin{eqnarray}
  \norm{x_{n+1} -x_n}^2 &= & \norm{x_{n+1} - T_n x_n}^2 - 2\psca{x_{n+1} - T_n x_n, x_n - T_n x_n} \nonumber \\
  & +& \norm{x_{n} - T_n x_n}^2  \nonumber \\
  & \ge &  \norm{x_{n} - T_n x_n}^2,
\end{eqnarray}
and we therefore obtain $\sum_{n \ge 0} \norm{x_{n} - T_n x_n}^2 < + \infty$. 
\end{myproof}


\end{document}